*Review*

# Inquiry-Based Mathematics Education: a call for reform in tertiary education seems unjustified

Tanya Evans[1*] and Heiko Dietrich[2]

[1]Department of Mathematics, University of Auckland, New Zealand; t.evans@auckland.ac.nz
[2]School of Mathematics, Monash University, Australia; heiko.dietrich@monash.edu
* **Correspondence:** t.evans@auckland.ac.nz; Tel: +64-9-923-8783

**Abstract:** In the last decade, major efforts have been made to promote inquiry-based mathematics learning at the tertiary level. The Inquiry-Based Mathematics Education (IBME) movement has gained strong momentum among some mathematicians, attracting substantial funding, including from some US government agencies. This resulted in the successful mobilization of regional consortia in many states, uniting over 800 mathematics education practitioners working to reform undergraduate education. Inquiry-based learning is characterized by the fundamental premise that learners should be allowed to learn 'new to them' mathematics without being taught. This progressive idea is based on the assumption that it is best to advance learners to the level of experts by engaging learners in mathematical practices similar to those of practicing mathematicians: creating new definitions, conjectures and proofs - that way learners are thought to develop 'deep mathematical understanding'.

However, concerted efforts to radically reform mathematics education must be systematically scrutinized in view of available evidence and theoretical advances in the learning sciences. To that end, this scoping review sought to consolidate the extant research literature from cognitive science and educational psychology, offering a critical commentary on the effectiveness of inquiry-based learning. Our analysis of research articles and books pertaining to the topic revealed that the call for a major reform by the IBME advocates is not justified. Specifically, the general claim that students would learn better (and acquire superior conceptual understanding) if they were not taught is not supported by evidence. Neither is the general claim about the merits of IBME for addressing equity issues in mathematics classrooms.

**Keywords:** explicit instruction; active learning; Inquiry-Based Learning (IBL); Inquiry-Based Mathematics Education (IBME); critical commentary

## 1. Introduction

Recent literature in higher education overwhelmingly depicts traditional lectures as an inferior method of instruction, counterposing it with a desired student-centered model in which active learning approaches are viewed as superior [1-4]. In the wake of the Covid-19 disruption to the sector, many tertiary institutions are considering mandating radical changes concerning instruction in response to repeated calls to abandon lectures altogether [1]. Many have already followed suit by replacing traditional lectures with flipped classrooms, given the meta-analytic evidence about the gained learning benefits [5]. The majority of evidence reported across numerous of studies points to the extra benefit of active engagement during class compared to passive listening, which is often experienced in transmission-style lectures.

However, under the banner of active learning, a new movement is gaining momentum in tertiary mathematics education – Inquiry-Based Mathematics Instruction [6]. In inquiry classrooms, students are challenged to "reinvent or create mathematics that is new to them" ([6], p. 131). A compelling premise of the approach is that in doing so learners engage in mathematics practices just like research mathematicians, through constructing definitions, formulating theorems and proving them. The underlying assumption is that in this process learners "develop deep mathematical understanding", coupled with a "sense of ownership through creation and reinvention" ([6], p. 131). The role of the lecturer is reduced to a facilitator of learning. The proponents of Inquiry-Based Mathematics Instruction consider inquiry a branch of active learning, thereby claiming the benefits associated with active learning reported in the research literature [1]. The similarity drawn with the active learning approach is justified on the premise that students in inquiry classrooms are expected to be actively engaged while reinventing mathematics. An inquiry-based approach, however, is characterized by the specific types of activities that require learners to engage in discovery on their own as opposed to following an explicitly explained method with provided worked examples. Inquiry curricula are structured with "a longer-term trajectory that sequences daily tasks to build toward big ideas" ([6], p. 131). An instructor designs task sequences to scaffold students' work on challenging problems so that students may prove a major theorem or (re)invent a mathematical definition or procedure over weeks of instruction.

The popularity of Inquiry-Based Mathematics Instruction has been growing, particularly in the USA. The Educational Advancement Foundation and its successor, Mathematics Learning by Inquiry, have played a foundational role in building the Inquiry-Based Learning (IBL) community and enabling evaluation research in the IBL Mathematics Centers [7]. A large grant from the National Science Foundation supported the growth of The Academy of Inquiry-Based Learning with a 5-year $2.8 million Collaborative Research Project, called Professional Development and Uptake through Collaborative Teams (PRODUCT), which started in 2015. The project funded a series of four-day intensive workshops bringing together over 300 mathematics faculty over the five-year period. A call for contributions for a special issue on 'Teaching inquiry' in Problems, Resources, and Issues in Mathematics Undergraduate Studies (PRIMUS) journal resulted in too many submissions to fit into a single issue, necessitating a second issue in order to publish 19 articles in total at the start of 2017 [8, 9]. A follow-up special issue, titled "Bringing Inquiry to the First Two Years of College Mathematics", was published in the same year (Issue 7), comprising a collection of articles promoting inquiry-based instruction in various settings, starting from pre-calculus [10]. In the same year, a Special Interest Group (SIG) on Inquiry-Based Learning was formed under the umbrella of the Mathematical Association of America (http://sigmaa.maa.org/ibl/). One of the aims of the group is "to promote the proliferation of IBL in Mathematics through conversation and professional development". The SIG has been efficient in providing mini-workshops and organizing symposia at conferences. For example, in 2018, at the Joint Mathematics Meetings of the American Mathematical Society, the SIG symposia on IBL comprised 50 talks over five half-days [6]. An establishment of purposeful fora to exchange ideas on IBL resulted in the annual National Inquiry-Based Learning and Teaching conferences, which started in 2018.

As a result of the coalescence of various initiatives supported by generous funding, the IBL community has grown into a large network of practitioners with centres based in different states across the USA [6]. Currently, a large project funded by the NSF provides targeted support for the development and support of the regional communities as part of the Communities for Mathematics Inquiry in Teaching Network (COMMIT Network) (https://www.comathinquiry.org/nsf-project). One of the goals of the project is to determine which initiatives are most effective in enabling sustained adoption of IBL in undergraduate mathematics education. So far, the project has successfully mobilized regional consortia in many states, uniting over 800 mathematics education practitioners working towards achieving a shared goal to improve mathematics education at the post-secondary level and enable success for all students.

In achieving this goal, however, concerted efforts to radically reform mathematics education must be systematically scrutinized in view of available evidence and theoretical advances in learning sciences. A nuanced, evidence-based approach to the umbrella term 'active learning' is required in order to discern the characteristics of learning activities that are more effective than others. To that end, the present study sought to consolidate the extant research literature from cognitive science and educational psychology, offering a critical commentary on the effectiveness of inquiry-based learning.

## 2. Critical Commentary on Inquiry-Based Learning

According to the key proponents, the goal of Inquiry-Based Learning is to transform students from consumers to producers of mathematics (as stated in the call for the MAA Contributed Paper Session on Inquiry-Based Learning and Teaching at the Joint Mathematics Meetings, 2020 https://www.jointmathematicsmeetings.org/meetings/national/jmm2020/2245_maacp-descrip). This progressive idea is based on the assumption that it is possible to advance learners to the level of experts via a shortcut of sorts by engaging learners in "mathematical practices similar to those of practicing mathematicians: conjecturing and proving, defining, creating and using algorithms, and modelling" ([6], p. 131). However seductive, this idea of discovery learning (or inquiry-based learning) has repeatedly been shown to be ineffective [11-17]. Despite the evidence available from numerous properly controlled, rigorously conducted experiments, many educational researchers remain strong advocates for the method due to the pivotal confusion – there is plenty of scientific evidence in support of 'active pedagogies' that are characterized by a learning environment in which students are attentively and actively engaged in their own learning. However, even though the so-called inquiry-based learning is considered a type of active learning, the scientific evidence on the benefits of active learning does not generalize to this very distinctive method. Until now, however, these approaches have advanced in the field of undergraduate mathematics education in parallel, with inadequate cross-referencing.

The current form of discovery/inquiry pedagogies can be traced back to the ideas of Jean-Jacque Rousseau, which reached us through many influential educators such as John Dewey, Jean Piaget, and Seymour Papert, to name a few [14]. For Rousseau and his successors, the unquestionable postulate is that it is always better to let the learners discover and construct the knowledge by themselves regardless of the time they spend tinkering and exploring. Rousseau optimistically believed that this time is never wasted. Quite the opposite; he assumed that this exploration would produce better-equipped minds capable of solving real problems rather than passively regurgitating received content knowledge to yield ready-made solutions acquired through rote learning [18]. "Teach your students to observe the

Writing:
4

phenomena of nature," says Rousseau, "and you will soon rouse his curiosity; but if you want his curiosity to grow, do not be in too great a hurry to satisfy it. Lay the problems before him and let him solve them himself", as cited in Dehaene ([14], p. 182).

The theoretical tenets of discovery/inquiry methods are rather compelling. However, a large body of research, accumulated over several decades, has demonstrated that its pedagogical value is unsatisfactory [13, 14]. The resounding failures of discovery learning replicated in so many studies prompted Richard Mayer, an American cognitive scientist, to write a review titled, "Should there be a three-strikes rule against pure discovery learning?" [19]. In this review, he analyzed the evidence from studies conducted from the 1950s to the late 1980s to compare pure discovery learning (unguided problem-based instruction) with guided forms of instruction. In conclusion, he made the following striking observation: in each decade since the 1950s, "after empirical studies provided solid evidence that the then-popular unguided approach did not work, a similar approach soon popped up under a different name with the cycle repeating itself" ([20], p. 6). Each new wave was driven by strong advocates who seemed unaware, at best, and dismissive, at worst, of previous evidence demonstrating that unguided approaches had not been validated. Such lamentable series of cycles produced discovery learning [21], which was superseded by experiential learning [22, 23], which gave way to problem-based and inquiry learning [24, 25], which were later joined by constructivist instructional techniques [26]. In summary, Mayer pointed out that the "debate about discovery has been replayed many times in education, but each time, the research evidence has favored a guided approach to learning" ([19], p. 18).

Given the mounting evidence of the ineffectiveness of discovery learning and its successive reincarnations [13], why is the approach still so prevalent and institutionalized in many educational settings, including initial teacher training colleges and university education departments, albeit only in those with a focus on sociological research? Historically, this is likely to be explained by the profound impact of one of the most influential and highly acclaimed cognitive psychologists, Jerome Bruner, who was undoubtedly swayed by the progressive ideas traceable back to Rousseau. In introducing the term 'discovery learning' in 1961 [27], Bruner wondered hypothetically about "the degree that one is able to approach learning as a task of discovery something rather than 'learning about' it" and whether there will be a tendency for the learner "to carry out his learning activities with the autonomy of self-reward or, more properly by reward that is discovery itself" ([28], p.17). Summarizing his thoughts, he stated that "the very attitudes and activities that characterize 'figuring out' or 'discovering' things for oneself also seems to have the effect of making material more readily accessible in memory" ([28], p. 24). This assumption was perfectly rational given the rudimentary state of development of cognitive psychology as a scientific discipline at the time. Therefore, it is unsurprising that the impact of Bruner's ideas was amplified through his scholarship as one of the founders of cognitive psychology. As a natural consequence, within two to three decades, this idea became dominant in educational spheres, seemingly unquestioned by the majority of education researchers. Furthermore, with the assurance of credibility rooted in cognitive psychology and the rapid growth of education faculty predominantly concerned with sociological aspects of education until recently, the situation persisted and remained largely unchallenged, especially in Anglosphere countries. In other words, one of the major reasons why so many educators mistakenly believe that minimally guided instruction is best for novice learners seems to be because they believe it is based on solid evidence from cognitive science.

However, as our understanding of human cognitive architecture progressed through the advances in learning sciences over the last six decades, it has become increasingly clear that discovery/inquiry learning is less effective than other methods [13]. Bruner's problematic theoretical justification was furthermore challenged by the accumulation of data from randomized, controlled trials, as well as from correlational studies strengthening the case against the inquiry learning perspective [13]. We next briefly outline the current understanding of the key aspects of human cognitive architecture relevant to instructional issues in the context of critical commentary on the merits of discovery/inquiry learning.

**2.1. Human cognitive architecture**

Our understanding of the mechanisms involved in learning has advanced substantially over the last few decades. Much research has been done to conceptualize those mechanisms by means of mathematical modelling and experiential testing of those models [29]. Diagrammatic representations of *Human Cognitive Architecture* [30, 31] can vary among researchers, but a modal model, which purposefully omits the complexity and simplifies the component relationships, can serve to explain the current state of research relevant to our context. A common way to represent human cognitive architecture is by using three memory systems: sensory memory, working memory and long-term memory. Figure 1, adapted from [32], depicts a modal model of human cognitive architecture, illustrating the main components and their interactions based on [29, 33, 34]. Sensory memory allows for the incoming sensory information (such as what we see, hear, touch, smell etc.) to be stored sufficiently long for the selected components to be transferred to working memory. Even after the stimuli have ceased, impressions of sensory information could be retained in sensory memory for short periods of time, provided they are selected for further processing in the working memory by the mechanism of attention [14].

Working memory is the domain where all conscious thoughts occur: if you are aware of something, it means this information is being processed in working memory. It has been known for some time that working memory has two critical limitations. First, the duration for which the information can be stored there is very short; almost all of the information is lost after 20 seconds without rehearsal [13, 35]. Thus, if you are working on novel information, you need to either write it down or keep rehearsing it until it is no longer needed. Moreover, there is now evidence that working memory depletes as a result of extensive mental effort before recovering after resting [36].

Second, working memory is extremely limited in capacity. It has been known since the mid-1950s that working memory can only hold about seven items of information (plus or minus two) [37]. However, only about three to four items of information can be processed at a time while mentally combining, comparing, or manipulating the items in some way [38]. As an example, on average, we can remember about seven random digits, but if asked to reorder them from, say, highest to lowest, the successful completion of the task would be challenging unless the number of digits is reduced.

In educational settings, the process of learning happens through working memory. There are two ways information can enter working memory. First, by focusing attention on particular incoming sensory information, one might consciously process a stimulus from sensory memory. According to Clark and Paivio's dual-coding theory [39] and Baddeley and Hitch's model of working memory [40], working memory has two systems: verbal/auditory processing takes place in the so-called *phonological loop*, whereas visual processing takes place within the *visuospatial sketchpad* (e.g.,[41]). The two systems are complementary, meaning that if processing can be split between the two systems, then the total working memory capacity could increase. Therefore, if a teacher presents instructional information by splitting it between visual and verbal modalities, then a learner's processing is more efficient (e.g., [34, 42]). This is, indeed, a customary instructional practice in mathematics teaching, where the explanations are usually presented by visual aids (such as written text with symbols or diagrams) with accompanying verbal narrations.

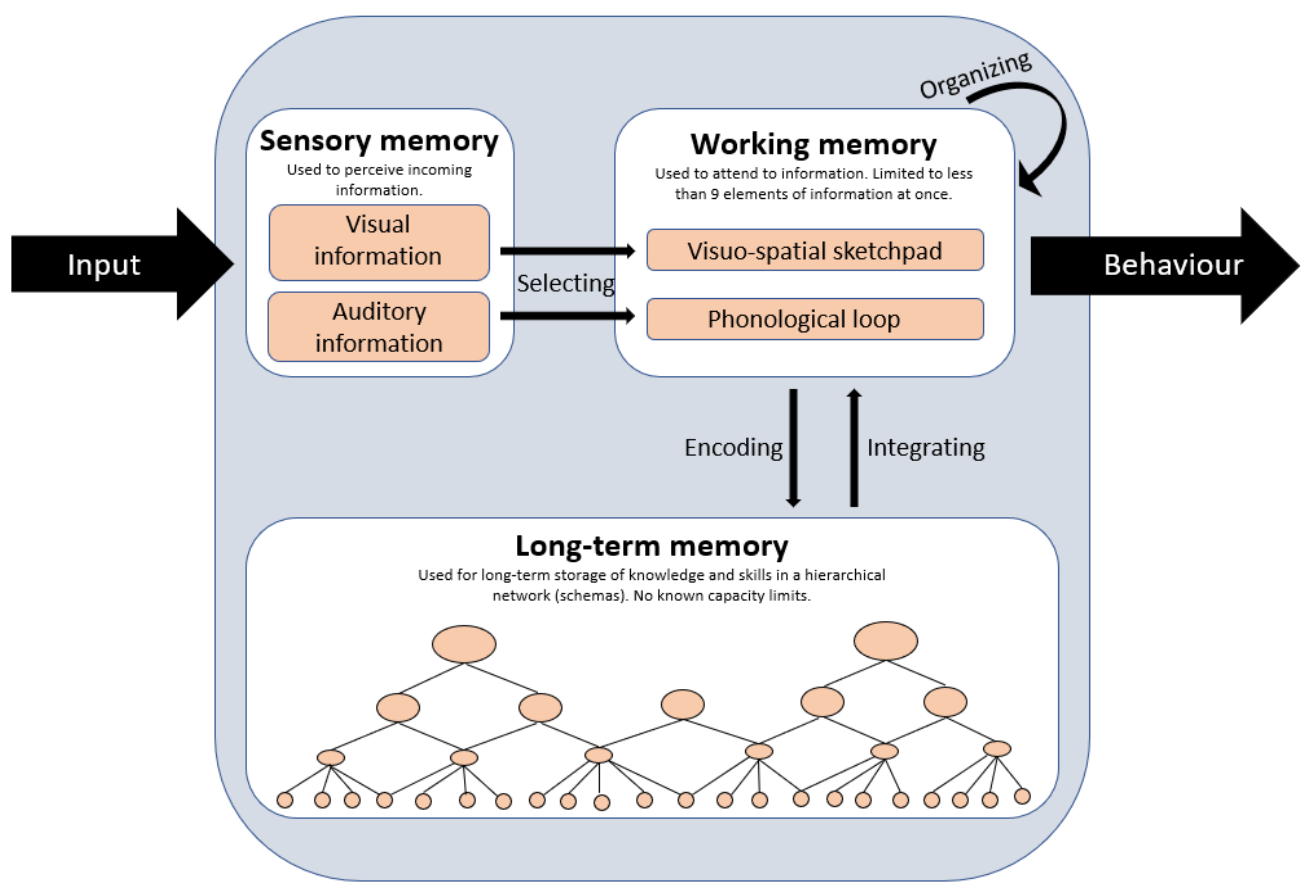

**Figure 1.** A modal model of Human Cognitive Architecture

The second way information can get into working memory is from the third memory system, long-term memory. Long-term memory is a central component in human cognitive architecture; it represents a repository of an enormous network of complex units of closely linked snippets of information, called schemas (Figure 1). We now know that our sense of self comes from that enormous amount of information stored in long-term memory [13]. Schemas are formed as outcomes of information processing in working memory. It is well-evidenced that experts have more comprehensive and better-organized schemas than novices, which results in improved fluency and accuracy in their performance. The surprising revelation that long-term memory plays a cardinal role in human cognition came from studies of chess players. It was long assumed that chess is a game of thinking and problem-solving. This turned out to be a fallacy: chess is a competition between the complexities of grandmasters' schemas stored in their long-term memory. Properly designed experiments, controlling for differences in working memory, showed that chess grandmasters have developed well-organized schemas that allow them to readily identify a myriad of board configurations [43, 44]. Those configurations are learned by studying previous games for many years. Chess grandmasters play well at a fast pace because all they are doing is recalling the winning configurations - not figuring out the moves [20, 44]. Similar results have been observed and validated in a wide range of educationally relevant areas from mathematics, computer science, physics, and social sciences [13, 45-49].

Hence, the main reason experts outperform novices is because of the superior quality of their schemas that organize a sizable body of knowledge and index it by a large number of patterns that, on recognition, guide the expert in a fraction of a second to the relevant parts [46]. These schemas constitute long-term memory and, when required, are retrieved and integrated into working memory [43]. As far as we know, at present, long-term memory appears to have no practical capacity limits and, thus, is used to overcome the temporal limitations of working memory. If you have stored a schema in long-term memory, you can then repeatedly reintegrate it into working memory, transcending the 20-second limit. Furthermore, long-term memory also functions as a bypass for dealing with the capacity limits of working memory. As stated before, only about three to four items of information can be processed in working memory at a time. However, if snippets of information are organized into a coherent schema and stored in long-term memory, then the whole schema can be brought into working memory as one item of information to be manipulated and integrated with other units of information. This way, the capacity of working memory is not exceeded, whereas it could be the case if the same knowledge had not been consolidated into a schema, for instance, if it is organized as disjoint schemas. For example, imagine someone tells you their cell phone number is 0149 1625 3649 and asks you to send them a message straight away. It could be challenging to hold this number in working memory while undertaking other tasks, such as pulling out your phone from your pocket, opening a messaging app etc. However, the same task would be straightforward if they told you that their phone number was a sequence of numbers squared from 0 to 7. In this case, the same information would have been organized into a coherent schema in your mind, thereby drastically reducing your working memory load [32].

**2.2. Inquiry-based learning from the perspective of cognitive science**

Advances in cognitive science have substantially altered our understanding of the role of long-term memory in human cognition over the last few decades. A large body of knowledge has developed since Atkinson and Shiffrin's [29] publication "Human memory: A proposed system and its control processes" in 1968, resulting in the development of an enormous research field, with 13,604 subsequent publications citing Atkinson and Shiffrin's work (retrieved 30 May 2022). However, the ideas and drivers of the movement that carried inquiry-based learning into the 21st century were predicated on the now-debunked assumptions about long-term memory. It is no longer viewed as a passive repository of discrete, isolated fragments of information that permit us to repeat what we have rote-learned [11]. And its influence on complex cognitive processes such as critical thinking and problem-solving is no longer viewed as subsidiary. Instead, long-term memory is now considered the central, principal structure of human cognition. Everything we see, hear, and think about is dependent on and influenced by our long-term memory. Numerous studies based on De Groot's (1945/1965) work on chess expertise, undertaken in the 1970-80s, have served as a major pivot for the reconceptualization of the field around the role of long-term memory.

The finding that expert chess players are markedly better than novices at reproducing briefly seen board configurations taken from real games, but do not differ in reproducing random board configurations, suggested that expert problem solvers derive their skill by drawing on extensive network of schemas stored in their long-term memory, which readily enable selection and application of the best procedures for solving problems [11].

The fact that the differences in the schemas in long-term memory could be used to fully explain problem-solving ability brought to the fore the importance of long-term memory in human cognition. We are skillful in something because our long-term memory contains a large network of organized schemas of snippets of information concerning the area. Those schemas enable us to quickly recognize the conditions of a situation and realize, often unconsciously, what to do and when to do it. Evolution has equipped us with the ability to store a vast amount of information in long-term memory, without which we would be struggling to complete even simple acts such as crossing a street. The schemas comprising long-term memory inform us how to avoid being hit by oncoming traffic - a skill many other animals are unable to encode in their long-term memories. Evolution is responsible for our ability to become competent in complex activities such as playing chess or doing mathematics through long-term practice and memorizing a countless number of 'problem states' and the best moves to make when encountering those states [11, 20, 49].

Since establishing the central role of long-term memory, educational psychologists and neuroscientists have devoted a great deal of attention to understanding how we learn most effectively. It is well-understood that in educational contexts learning is not just an additive process of storying new snippets of information into memory, as in a computer [33]. Rather, learning is dependent on two major factors: (1) what is presented and how and (2) the cognitive processing that the learner is actively engaged in during learning. Thus, learning is viewed as a generative activity. This conception of learning is a natural theoretical advance built on the coalescence of cognitive revolution and constructivist ideas about the importance of meaning-making while learning. This well-evidenced theory envisions effective learning to comprise three stages: (1) learners actively select the relevant aspects of incoming information by paying attention, (2) which is followed by organizing this information into a coherent cognitive structure in working memory, and (3) integrating cognitive structures with relevant prior knowledge activated from long-term memory (see Figure 1) [33]. Hence, the effectiveness of a learner's cognitive processing plays a pivotal role in generative learning. This does not mean that effective learning must necessarily involve a behaviorally active component. There is a functional difference between behaviorally active engagement and cognitively active engagement – these are two different variables. In fact, research shows that behaviorally active engagement does not necessarily lead to effective learning [19], whereas behaviorally passive engagement, such as viewing a power-point presentation, can lead to active learning, provided the passive instruction is designed well (reducing extraneous processing, managing essential processing, and fostering generative processing [50]). The conflation of the functional roles of these two types of active engagement in the learning process is known as the "constructivist teaching fallacy", which occurs when a teacher assumes that active learning is caused by active instructional methods and passive learning is caused by passive methods of instruction [19].

It appears that the proponents of pedagogical approaches under the umbrella of inquiry-based learning have fallen for the "constructivist teaching fallacy". Or they may only be aware of the outdated model of the human mind, according to which minimal guidance problem-solving through inquiry is the most effective way to gain expertise. In other words, it appears that the advocates for inquiry-based learning miss to account for the role of long-term memory in learning and also neglect to employ a nuanced consideration of active (cognitive) engagement involved in learning processes. This is evidenced by their commonly used negative connotations in arguing against the "traditional, fact-oriented, teacher-centered model" of instruction [50] and erecting a false dichotomy between learning for knowing versus learning for conceptual understanding. We will return to this point in the next section.

First, we need to clarify the terms. By inquiry-based learning in undergraduate mathematics we mean a specific approach when instructors expect students to discover 'new to them' mathematics by asking them to formulate definitions and theorems, identify algorithms and prove theorems. Instructors act as facilitators by providing an organized instructional sequence of tasks that students are expected to complete during class (sometimes over many sessions), usually as part of a group. Group discussions are strongly encouraged [6]. However, we want to distinguish this from assigning a non-routine task (for which no analogous worked example is provided) at the end of a topic to facilitate student inquiry to consolidate the learned material that was adequately explained (or modelled) by an instructor. In fact, we contend this would be a beneficial activity, according to the generative account of learning. This is because it may well be the case that learners would have developed relevant schemas in their long-term memory as a result of teacher explanations, worked examples and solving/proving/working on similar problems assigned to them, but their schemas might still be disjoint. Such inquiry tasks could successfully facilitate the unification of those schemas into one high-quality schema, thus improving the organization of structure in the long-term memory. This could result in a positive change in long-term memory, reflecting improved understanding. However, we will argue that this could only be effective if a learner has already acquired sufficient schemas in their long-term memory about the topic.

Next, we will present a critical commentary pertaining to the three major claims often made by inquiry-learning proponents that: (1) students learn mathematics better when they are not taught (explanations of new material should be minimized or absent; students are to engage in doing/discovering mathematics by themselves mimicking mathematicians); (2) inquiry-based learning results in quality learning outcomes (superior conceptual understanding, which enables the transfer of learning to new situations), and (3) inquiry-based classrooms are more beneficial for addressing equity issues than conventional classrooms.

**2.3. Novice students do not learn better when they are not taught**

The overarching sentiment that teaching is not an act of goodness seemed to be rooted in the amalgamation of the enduring traditions of the progressive education movement based on the ideas traceable back to Rousseau and the ultra-progressive sociopolitical ideas that started to dominate educational scholarship at the turn of the 20th century. An act of teaching has been construed as an expression of repressive power imbalance. For example, in a publication by the Australian Association for Research in Education, teaching (direct instruction) has been criticized because "it places the teacher and child in a rigid relationship where the teacher is always the one with the power and knowledge with limited allowance or recognition of individual and cultural difference" ([51], p. 2).

Very often, the righteousness of this proclamation has been asserted through invoking value statements introducing research such as "This paper sought to look at the mathematics teachers' effort to shift from the traditional teacher-centered classroom instruction to a democratic student-centred classroom." ([52], p. 1). The implication is that a teacher-centered classroom is simply undemocratic and does not align with our value system.

It appears that in a perplexing twist of confusion between different meanings of the word "progressive", the idea that a teacher should be a non-authoritative facilitator of learning has been propelled by well-meaning, empathetic educators on the premise of progressive political agendas. However, as has been pointed out by Greg Ashman in his book [53], progressive education is not progressive politics. Rather, progressive education is a teaching philosophy with the core assumption that learners need to be accommodated in a naturalistic way of learning through following their interests. This philosophical position is not inconsistent with any point on the political spectrum. Interestingly, noting that many left-leaning educators feel obliged to adopt such progressive ideas as default, progressive education was once the chosen philosophy of Giovanni Gentile, Mussolini's education minister. Unsurprisingly, as a natural counterpoint, it drew fierce criticism from the Marxist philosopher Antonio Gramsci [53, 54], underscoring the individualistic focus of progressive education, which is usually associated with the political right. Indeed, things tend to get very convoluted when political ideologies get entangled with philosophical positions in education, which, in turn, manifest in reform agendas in the classrooms. Unfortunately, it appears that this confusion deeply permeated academic and bureaucratic educational circles in the last few decades, particularly in the Anglosphere. The idea of a teacher providing explicit instruction to learners is a deeply offensive, regressive idea to many [53]. Perhaps, such a pervasive stance afforded justification for many not to question or seek evidence on the effectiveness of the progressive pedagogical approaches.

Within the education research community, the issue has come to a head with the influential publication by Kirschner, Sweller and Clark, titled "Why Minimal Guidance During Instruction Does Not Work: An Analysis of the Failure of Constructivist, Discovery, Problem-Based, Experiential, and Inquiry-Based Teaching" in 2006. They presented a strong case against discovery methods of instruction by summarizing "empirical studies over the past half-century that consistently indicate that minimally guided instruction is less effective and less efficient than instructional approaches that place a strong emphasis on guidance of the student learning process" ([11]p. 75). In the following year, a live debate was held at the 2007 annual convention of the American Educational Research Association, which continued in a book, "Constructivism Instruction: Success or Failure?" edited by the organizers of the debate, Sigmund Tobias and Thomas Duffy [50]. Many prominent scholars on both sides of the debate of the controversy on the success or failure of constructivist instruction contributed a chapter to this volume. After completing their chapters, authors on both sides responded to questions from two scholars on the other side of the debate, with several iterations of questions and responses in some cases. This 376-page volume presented a comprehensive overview of the existing research, evidence, and arguments presented for both sides of the debate. In conclusion, Sigmund Tobias stated: "A careful reading and re-reading of all the chapters in this book, and the related literature, has indicated to me that there is stimulating rhetoric for the constructivist position, but relatively little research supporting it. Indeed, the constructivists seem to suffer from denial with regard to information processing" ([50], p. 346). Specifically for mathematics education, he stated: "In my assessment, the preponderance of evidence shows the ineffectiveness of minimally guided methods for promoting mathematical understanding" ([50] p. 285).

This conclusion was based solely on the available evidence to date and the explanation afforded by the information processing theories based on Human Cognitive Architecture. In order to reap benefits from inquiry-based learning, a learner needs to select, organize, and integrate high-level information in a task-appropriate way, which is quite demanding for novice learners. Research has identified that inquiry/discovery learning relies on an extensive search through problem-solving space and that such a process is likely to tax learners' limited working-memory capacity and thus is unlikely to result in efficient learning for novices [55, 56]. An extensive body of expert-novice research has shown that novices cannot recognize what is relevant to a problem or problem solution and also struggle to recognize what is novel in a specific situation. In other words, what you know determines what you recognize [50]. Moreover, research has shown that, left to their own devices, novices tend to focus on surface features instead of structural features of problems when solving problems and selecting further problems for study and solution [57]. In summary, according to cognitive load theory, the exploration of complex phenomena by novices with minimal guidance is likely to impose heavy loads on working memory, thus resulting in detrimental learning outcomes [11, 13, 55, 58-60].

Furthermore, in inquiry-based learning of novel information, learners would need to efficiently monitor and control their own processes of attention to relevant information [11, 14, 61]. This would require learners to have well-developed metacognitive skills, which is an unwarranted assumption to be made about all learners [62]. In summary, from a theoretical perspective, learning by inquiry/discovery would necessitate a greater number of mental operations and require better executive control of attention compared to learning under a more directive approach (when an expert explanation is provided). For example, suppose a lecturer offers an explanation using both auditory/verbal and visual channels to reduce the processing load of the working memory (referrer to Figure 1). In this case, the explanation would help focus learners' attention on the incoming information so that the learner can be actively (cognitively) engaged in selecting the relevant information for further processing in working memory. In other words, learners are enabled to process this information in the working memory actively and integrate it with the existing schemas activated from the long-term memory to be re-encoded in a modified, extended way. Thus, a provision of an expert explanation is conducive to an efficient way of learning novel information.

This explanation gained further empirical support when two meta-analytic studies were published in 2011 by Alfieri et al. [62]. The first one meta-analyzed the findings of the existing studies that directly tested for differences between an explicit instruction condition and a condition in which unassisted discovery learning occurred. Unassisted discovery learning condition was operationally defined as being provided with no guidance or feedback during the learning task (also known as pure discovery/inquiry learning). The inclusion criterion was strict only to include articles reporting on comparable conditions that consistently differed in the type of instruction. Studies comparing fundamentally different or equivocated conditions prior to testing were not included. Qualitative studies were also excluded from the meta-analysis due to the limitations in providing generalizable results.

In total, 580 comparisons from 108 studies were included in the meta-analysis (some studies reported multiple experiments). Utilizing the random effects analysis, the 108 studies had a mean effect size of $d = -0.38$ (95% CI [–.50, –.25]), demonstrating that explicit teaching was more beneficial to learning than unassisted discovery (with small but meaningful effect size, $p < .001$). Moderation analysis indicated that the type of publication moderated the findings, with the articles in first-tier journals ($d = –0.67$) reporting larger effect sizes in favor of explicit instruction than articles in second-tier publications ($d = –0.24$). Effect sizes from book chapters were the smallest reported: $d = –0.12$. Furthermore, participants in unassisted discovery fared worse than participants in all four comparison conditions comprising explicit instruction: direct teaching ($d = –0.29$), feedback ($d = –0.46$), worked examples ($d = –0.63$), and explanations provided ($d = –0.28$). Therefore, the findings indicated that explicit-instructional conditions lead to greater learning than unassisted-discovery conditions. Moreover, there were no significant differences between the mean effect sizes of the three categories of unassisted-discovery conditions: unassisted tasks, tasks requiring invention, and tasks involving collaboration with a naive peer were all found to be equally detrimental to learning.

The second meta-analysis by the same authors [62], employed a different perspective with the focus of comparison being the enhanced discovery approach which comprised three categories: elicited explanations (e.g., self-explanation prompts), guided discovery (some instruction) and generation (similar to the invention), which were compared to other instructional methods (covering the entire spectrum from unassisted discovery to direct teaching). A total of 360 comparisons from 56 studies indicated that the type of enhanced-discovery condition moderated the findings: elicited explanation ($d = 0.36$) and guided discovery ($d = 0.50$) favored enhanced discovery, whereas generation ($d = –0.15$) favored other instructional methods. The latter finding was unexpected given the known benefits reported as the generation effect [63]. For instance, given pairs of synonyms, experiment participants will remember a word better if they explicitly have to generate missing information, as in the case of "FAST: R_P_D" versus reading "FAST: RAPID". The generation effect is often used as a kernel of an argument in the discovery/inquiry literature. However, Alfieri et al. found that generation is not an optimal form of enhanced discovery. "The generation conditions required learners to generate rules, strategies, or images or to answer questions about the information" ([62], p. 12) – tasks largely similar to the activities used in the Inquiry-Based Mathematics Instruction classrooms [6].

It is fair to say that this meta-analytic study has put an end to a long debate. The evidence is that the unassisted form of discovery is inferior to any other type of instruction and that the addition of scaffolding and guidance (such as expert explanations and self-explanation prompts) leads to better outcomes. The big idea of discovery pedagogy nurtured by the followers of Rousseau and embellished during the 20[th] century by the virtuous progressive educationalists has been decisively discredited by the mounting experimental evidence. In their defense, the idea that learners should be able to construct their own understandings of academic disciplines with minimal assistance because they do so on a daily basis in the context of everyday activities is highly compelling. And if it turned out to be true, it would have most definitely revolutionized education and optimized learning for all involved. That said, it made perfect sense for the idea to be properly developed and put to the test. Unfortunately, it did not work. A possible explanation why is offered by the theory of evolutional psychology [64], which draws the distinction between evolutionary-primary knowledge (how to walk, talk, eat, and solve problems) that we acquire almost effortlessly by being emerged in a social environment and secondary (academic) knowledge that we have not evolved to acquire so easily. From the evolutionarily perspective, the content and context of formal education are extraordinary [65] and thus require more explicit assistance to arrive at accurate knowledge constructions, understandings, and solutions [12].

Given the evidence from numerous properly controlled experiments and a theoretical explanation of the inferiority of the discovery pedagogy, a large majority of educational researchers who are familiar with educational psychology have moved on. However, the remaining proponents of the pure discovery/inquiry method seemed to be unaware of the past existence and the closure of one of the major educational debates. Currently, the main open questions debated pertain to the type of explicit assistance (e.g., quality of explanations) [66-69], scaffolds (e.g., self-explanation prompts) [70-74], and their optimal sequence in combination with practice problems to enable effective learning [75].

In terms of explicit assistance, the effectiveness of instructional explanations is now back in the research foci, with some studies attempting to conceptualize their explanatoriness [68] and identify features making them most effective [66]. Comparing the quality of explanations provided by mathematicians (with lower pedagogical content knowledge but high content knowledge) and mathematics teachers (with high pedagogical content knowledge but lower content knowledge) about an extremum problem intended for high school students, Lachner and Nückles [76] concluded that the explanations were principally different. The teachers mainly described the solution steps as an algorithm for finding extreme values of a function (product-orientation). However, in addition to demonstrating the solution steps, the mathematicians provided clarifying information explaining why a certain step in the solution was necessary (process-orientation). In the follow-up study, the researchers investigated the effectiveness of these explanations on eighty high-school students who were randomly split into three groups, receiving: (1) product-oriented explanations, (2) process-oriented explanations, or (3) no explanations (discovery group). Consistent with other research, they found that students who were not provided an explanation showed the lowest learning gains. Moreover, students who learned with a process-oriented explanation performed significantly better than students who learned with a product-oriented explanation on an application test. The latter finding was replicated and extended in a 2019 study [67] in which the impact of different explanations was assessed in a randomized sample of 129 students receiving either principle-oriented or procedure-oriented explanations on four mathematical topics. Students who were given principle-oriented explanations substantially outperformed students given procedure-oriented explanations on the application test (with similar problems) and the transfer test (with dissimilar problems).

Overall, the evidence from experimental psychology research does not support one of the main claims made by the Inquiry-Based Mathematics Education proponents that students learn better when they are not explicitly taught. The benefit of being taught by an expert is explained by cognitive science since experts' well-integrated and highly sophisticated knowledge [14, 77, 78] helps to produce explanations that can serve novice students as powerful cognitive scaffolds for their thinking and learning processes [79]. Research on expertise showed that experts mainly organize their knowledge around abstract principles of a domain, whereas novices' knowledge structures are often lacking integrative organization [76-78]. Therefore, a globally cohesive expert explanation, which focuses on principles, helps novices to mentally organize the presented information into a coherent schema, thus enabling them to integrate the information into their prior knowledge. As long as a learner is paying attention to the suitable explanation [14], it is predicted that the learner would be actively (cognitively) engaged in generative processing, resulting in effective and efficient learning.

**2.4. Inquiry-based learning is not the most effective way to generate conceptual understanding**

Advocates of inquiry-based learning concur with Piaget's assertion that "each time one prematurely teaches a child something he could have discovered for himself, that child is kept from inventing it and consequently from understanding it completely" ([80] p. 715). Moreover, they argue that those who acquire knowledge without explicit instruction are more likely to apply and extend that knowledge [81, 82]. Proponents of Inquiry-Based Mathematics Education (IBME) at the tertiary level assert that instead of receiving an explanation from an expert, students should be engaged in tasks similar to those of practicing mathematicians: defining, conjecturing, proving, creating and using algorithms, and modeling. In doing so, IBME advocates claim, "students not only develop deep mathematical understanding, but they also develop a sense of ownership through creation and reinvention" ([6], p 131). The central theme is that inquiry-based learning is more effective for developing deep conceptual understanding than explicit instruction.

What is conceptual understanding? An understanding of a mathematical concept is formed when a learner has encoded a large schema relating the concept with other preexisting schemas of related concepts and their relationships. The bigger the schema, the better the understanding. Experts have large, highly-integrated, and organized schemas connecting many subschemas in one conceptual interpretation, which could be referred to as 'understanding'. Efforts to conceptualize and measure conceptual understanding comprise an active area of research in mathematics education. A major challenge is to measure knowledge of a given concept with acceptable validity and reliability [83]. However, a widely accepted proxy for measuring conceptual understanding is to examine a learner's performance on transfer problems – a set of non-routine questions about a concept/procedure that were neither previously presented as worked examples nor assigned as practice exercises.

A pivotal study examining the effect of different instructional methods on transfer was conducted by Klahr and Nigam [84] in the context of the control-of-variables strategy (CVS). To master CVS, elementary school children needed to learn how to create experiments in which a single contrast is made between experimental conditions and logically understand the inherent indeterminacy of confounded experiments. In other words, acquiring CVS is a crucial step in developing scientific reasoning, enabling children to design unconfounded experiments from which valid causal inferences can be drawn. In their experimental study (N = 112) with random assignment to two groups, the researchers replicated previous studies showing that direct instruction was clearly superior to discovery learning in facilitating children's acquisition of CVS [85, 86]. Direct instruction resulted in significantly more masters of CVS than discovery learning (77% in the direct-instruction condition vs 23% in the discovery learning condition). Importantly, the superiority of direct instruction remained after restricting the analysis to only low-performing students, with 69% of the 35 in the direct-instruction condition becoming masters, compared with only 15% of the 41 children in the discovery condition with equally low initial scores. However, the most important result of the study demonstrated that the learners who became masters via direct instruction were as skilled at the transfer task (evaluating science-fair posters) as were discovery-learning masters and experts. Additionally, learners who did not become masters did equally poorly on the poster-evaluation transfer task regardless of the instruction condition. Hence, the overall result is rather prominent: the focused, explicit, and didactic training in the direct-instruction condition resulted in a large majority (77%) of learners becoming CVS masters who were as proficient as relatively few (23%) discovery-learning masters when subsequently asked to demonstrate conceptual understanding in an authentic environment, offering scientific judgments. The researchers concluded with a call "to reexamine the long-standing claim that the limitations of direct instruction, as well as the advantages of discovery methods, will invariably manifest themselves in tasks requiring broad transfer to authentic contexts (e.g., "learning under external reinforcement… produces either very little change in logical thinking or a striking momentary change with no real comprehension''—Piaget, 1970, p. 714)" ([84], p. 666).

The effect on transfer in mathematics contexts was investigated by Lachner, Weinhuber, and Nückles in the 2019 study, mentioned earlier [67]. In a randomized sample of 129 students learning to solve extreme value real-world problems, students who were given principle-oriented explanations substantially outperformed students given procedure-oriented explanations on the transfer test. The researcher offered a theoretical explanation based on the findings of related studies about comparing expert and advanced-student explanations and their impact on novice medical students [78]. They found that the experts' explanations were more integrated than advanced-student explanations. Qualitative analysis of the explanations' constituent elements (concepts) and their connections to others (similar to concept-mapping) revealed that experts' elements were more thoroughly interconnected, yet they contained fewer details than the advanced students' explanations. They also found that the experts' explanations resulted in superior novices' learning outcomes when the novices were given transfer tasks that required a reinterpretation rather than a straightforward application of the previously acquired knowledge. In the next study, using a detailed analysis of how students engaged with the explanations, they observed that novices studying an expert's explanation tended to follow a deep-processing approach, whereas novices with an advanced student's explanations processed the explanations in a shallow way. Novices studying an expert's explanation were enacting self-explanation loops more frequently. In contrast, this deep-processing sequence was almost completely absent in the group of novices learning with an advanced student's explanation. The conclusion of this research is that the extent to which explanations were conducive to novices' acquisition of transferable knowledge is heavily dependent on the instructor's level of domain expertise. This is because experts' explanations tend to be globally cohesive, directing learners' attention to key concepts and principles, thus enabling the realization of a deep-processing approach to learning.

On the other side of the debate, there are a number of publications alluding to benefits afforded by discovery/inquiry-based approaches that stress the potential positive impact on the transfer of learning. However, at the time of the 2007 AERA debate and its subsequent continuation in the book "Constructivism instruction: success of failure?", none of them were deemed to provide a conclusive evidence base because those studies did not employ suitable methodologies to allow for such conclusions [50].

Since then, the value of explicit instruction for enabling efficient learning of novel information as part of a learning episode has been widely accepted. As a result, the research debate has shifted to questions about optimal sequencing within a learning episode. At the top of current debates is an important question: whether or not a problem-solving activity (by learners) is beneficial prior to an explicit explanation provided by an expert. An approach termed 'productive failure' has become popular with practitioners and researchers in which learners first struggle to solve a problem on their own before being provided with an explicit explanation on how to solve the problem [87]. Kapur [88] suggested that problem-solving first may activate and differentiate learners' prior knowledge, thereby increasing their awareness of the gaps in their prior knowledge. That way, learners might be better able to attend to the critical features of the explanation provided. Moreover, Kapur suggested that learners involved in problem-solving first may be more motivated and engaged. Additionally, other researchers hypothesized that "requiring learners to generate their own problem solutions prior to explicit guidance may strengthen the stimulus-response relation in memory in a similar way as has been proposed in order to account for the 'generation effect' [63, 89, 90]" ([75], p. 230). To test this plausible hypothesis,

Ashman, Kalyuga and Sweller [75] conducted a fully randomized, controlled experiment in which learners were randomly assigned to one of two conditions: (1) a problem-solving–lecture sequence and (2) a lecture–problem-solving sequence. They found that learners in the lecture–problem-solving condition (explicit instruction first) scored significantly higher not only on an application test but also on a transfer test than learners in the problem-solving-lecture condition (problem-solving first). Mean test scores were almost 50% higher for the lecture–problem-solving sequence for both tests. Hence, they concluded, for learning where element interactivity is high, explicit instruction should precede problem-solving.

Although some evidence base exists in support of 'productive failure', Ashman et al. [75] identified *element interactivity* (complexity of a concept in terms of its connections to and dependence on other concepts) as an explanatory variable for the presence or absence of a problem-solving first advantage. This is in line with the expertise reversal effect, which is now considered a variant of the more general element interactivity effect [91]. That is, to learn high element interactivity information, studying worked examples is more beneficial than problem-solving. However, with the increase in expertise, the element interactivity of the information decreases (as an outcome of learning); hence the advantage of worked examples reduces and eventually reverses [92, 93]. In a tertiary mathematics context, most novel concepts and procedures are arguably information with a high level of element interactivity. Therefore, it follows that in order to support conceptual understanding and transfer, a learning episode should begin with explicit teaching (providing explanations) before shifting to more problem-based instructional methods.

## 2.5. Inquiry-based classrooms are not more beneficial in addressing equity issues than conventional classrooms

In line with the progressive education ideas, inquiry-based mathematics classrooms have been portrayed as a panacea for all equity issues [94]. The main assumption is that inquiry classrooms are equitable environments where educators can "face specific concerns about whose voices are privileged or excluded in mathematics" and "recognize similar issues of identity, agency and power in their own higher education settings" ([6], p. 140). The concept of 'personal empowerment' as a form of transformative learning has been introduced in mathematics education literature [95]. It is meant to encapsulate the enhancing elements of inquiry-based mathematics learning through consideration of self-empowerment, cognitive empowerment, and social empowerment. Using semi-structured interviews with undergraduate students, researchers argue that participation in inquiry-based learning could be strongly transformative for individual students. They assert that "not only do these courses enhance students' thinking and problem-solving skills but they also significantly promote self-perceptions, agency and self-regulatory activity, and social skills" ([95], p. 316). However, evidence from properly controlled randomized experiments to support such claims is missing.

Within undergraduate mathematics education literature on inquiry-based learning, one of the most cited studies by Laursen et al., [96] published in one of the top journals, the Journal for Research in Mathematics Education (JRME), in 2014, claimed that inquiry-based classes level the playing field for women by offering learning experiences of equal benefit to men and women. This conclusion was reached based on the analysis of self-reported learning gains reported by students from 42 mathematics sections taught in inquiry-based format, which were compared to conventional courses taught at four US universities with established centers for Inquiry-Based Learning. Apart from the apparent limitation of relying on self-reported learning gains data, the researcher did not compare learning outcomes by altering only one variable at a time (to control for confounding variables). Rather, inquiry-based sections were taught by enthusiastic educators with a buy-in for the cause, which would have unavoidably affected the outcomes. Moreover, a serious methodological flaw was that participants were not randomly assigned to different conditions, which could have resulted in self-selection with biased characteristics. Despite the limitations, the authors concluded that non-inquiry-based courses "do selective disservice to women" ([96], p. 415).

This problematic conclusion was exposed by a series of recent studies which contradicted Laursen et al.'s [96] claims. Johnson et al. (2020) [97] examined the relationship between gender and student learning outcomes in inquiry-oriented abstract algebra. Using hierarchical linear modeling, content assessment data were analyzed from 522 students, which identified a gender performance difference. Men outperformed women in the inquiry-oriented classes, whereas no such gender difference was detected in conventional classes. This finding was corroborated by the latest study by Reinholz et al. (2022), which was also published in JRME, titled "When Active Learning Is Inequitable: Women's Participation Predicts Gender Inequities in Mathematics" [98]. The authors state that contrary to the conventional wisdom that inquiry-based learning promotes equity, "across a sample of 20 undergraduate mathematics classrooms, we found evidence of greater gender inequity in favor of men in the inquiry-oriented instructional environments" ([98], p. 218). In contrast, no significant performance differences were found between genders in the non-inquiry comparison samples. Using a weighted regression analysis, the researchers explained the findings by identifying a significant link between women's participation rates during classroom discussions and gendered performance differences. The findings uncovered that women tend to be more affected by social aspects of learning than men, which, in turn, translates into tangible effects on learning outcomes [98].

Overall, from a cognitive theory perspective, inquiry-based mathematics classrooms are predicted to be less effective for learners with low prior knowledge, which is often a characteristic of learners from underrepresented groups. As mentioned previously in the Human Cognitive Architecture section, when processing novel information (and only when processing novel information), working memory has two major constraints: it is extremely limited in both capacity and duration. However, these limitations disappear if the information is accessed from long-term memory [31]. Taken together, these two facts provide a compelling argument explaining why partially or minimally guided instruction is generally ineffective for novices, despite often being effective for experts. In other words, novices only resource when given a mathematical novel task is their seriously constrained working memory [20]. In contrast, experts have both their working memory and all the relevant knowledge and skills stored in long-term memory, which are readily accessible for processing in working memory.

Hence, inquiry learning or any instructional activity that unnecessarily increases working memory load will inevitably have deleterious consequences [11]. This theoretical explanation is supported by a large number of properly controlled experiments reported in the research literature in the last 60 years [11, 13, 55, 58-60, 99]. In summary, cognitive theory and experimental empirical evidence present a serious challenge to the claim often made by advocates of inquiry-based learning that minimal guidance instruction is beneficial for addressing equity issues.

## 2.6. Conclusion

In this commentary, we leveled our criticism at a particular instructional practice, often termed

inquiry- or discovery-based learning, which is characterized by the fundamental premise that learners should be allowed to learn 'new to them' mathematics without being taught [6]. This manifests in a change of the lecturer's role to be a 'facilitator of learning' instead. The role of an inquiry-oriented instructor is to set up inquiry tasks that learners are expected to complete in order to discover new-to-them knowledge, often as a part of small group activity. In the last decade, major efforts have been devoted to promoting the adoption of inquiry-based mathematics learning at an undergraduate level. The movement, Inquiry-Based Mathematics Education (IBME), has gained momentum with hundreds of participating mathematicians and attracted substantial funding from the US government agencies.

However, given the preponderance of evidence summarized in this commentary from experimental research and the theoretical accounts explaining how efficiency in learning is best achieved, we conclude that the call for a major reform from the IBME advocates is not justified. Specifically, the general claim that students would learn better (and acquire superior conceptual understanding) if they were not taught is not supported by evidence. Neither is the general claim about the merits of IBME for addressing equity issues in mathematics classrooms. In summary, the claim that students should be allowed to discover mathematics for themselves akin to professional mathematicians conflicts with a large body of research supported by experimental evidence. Therefore, inquiry-based education should not be promoted for general adoption in mathematics undergraduate education without further research. Our recommendation, instead, is to consider ways to promote active (cognitive) engagement of learners in a conventional mathematics lecture. The provision of explanations by an expert should be supplemented by prompts for students to engage in self-explanations during a learning episode (be it a face-to-face or online lecture). This could be achieved in many different ways. For example, quizzes with self-explanation prompts, targeted questions from a lecturer to promote self-explanations, or short small group tasks could be used as part of or in the conclusion of an instructional explanation during a lecturing episode. Research has shown that such activities are likely to engage learners in generative learning and are intended to prime appropriate cognitive processing during a lecture, such as paying attention to the relevant information, mentally organizing it, and then integrating it with the learners' prior knowledge [14, 33]. Future research should investigate various ways to improve conventional mathematics lectures (face-to-face or online videos) by incorporating prompts for active (cognitive) engagement of learners as part of the explanation sequences.

**Author's biography**

**Dr. Tanya Evans** is a senior lecturer in the Department of Mathematics at the University of Auckland. She is specialized in mathematics education. Her research interests include secondary and undergraduate mathematics education, mathematical practice, professional development, and curriculum studies.

**Dr. Heiko Dietrich** works at Monash University as an associate professor in the School of Mathematics and the Associate Dean Graduate Education in the Faculty of Science. His research interests include computational algebra and aspects of algebraic design theory, in particular, finite groups, groups of prime-power order, and Lie algebras.